\documentclass%[amsmath]
{article}
\usepackage{latexsym}
\usepackage{amsthm}
\usepackage{amsfonts}
\usepackage{amssymb}

\theoremstyle{plain}
  \newtheorem{theorem}{Theorem}[section]

  \newtheorem{lemma}{Lemma}[section]

\theoremstyle{remark}
  \newtheorem{remark}{Remark}[section]

\theoremstyle{definition}
  \newtheorem{definition}{Definition}[section]
  \newtheorem{notation}{Notation}[section]

\makeatletter

\@addtoreset{equation}{section}
\makeatother

\begin{document}

\title{
Local well-posedness of the nonlinear Schr\"odinger equations on the sphere for data in modulation spaces}
\author{Hideo Takaoka\thanks{Department of Mathematics, Hokkaido University, Sapporo, 060-0810, Japan. 
Emai: takaoka@math.sci.hokudai.ac.jp
The author was supported in part by J.S.P.S. Grant-in-Aid 25287022.}}

%\pagestyle{myheadings}
%\markboth{H. Takaoka}{}

\date{}

\maketitle

\begin{abstract}
In this paper we discuss a priori estimates derived from the energy method to the initial value problem for the cubic nonlinear Schr\"odinger on the sphere $S^2$.
Exploring suitable a priori estimates, we prove the existence of solution for data whose regularity is $s=1/4$.
% and the (unconditionally) local-in-time well-posedness in Sobolev spaces $H^s$ with $s>3/4$.
\end{abstract}

\section{Introduction}

In this paper we consider the initial value problem for the cubic nonlinear Schr\"odinger equation on the surface:
\begin{eqnarray}\label{NLS}
\left\{
\begin{array}{l}
i\partial_t u+\Delta u=|u|^2,\\
u(x,0)=u_0(x),
\end{array}
\right.
\end{eqnarray}
where the unknown function $u=u(x,t)$ is a complex valued function on $(t,x)\in\mathbb{R}\times S^2$, and $\Delta$ denotes the Laplacian on the surface of the unit sphere $S^2$ in $\mathbb{R}^3$.
The equation (\ref{NLS}) is defocusing.
However the result in this paper is independent on whether the equation is focusing or defocussing.

In the case of the whole space case $\mathbb{R}^2$, local well-posedness in $H^s(\mathbb{R}^2)$ for $s\ge 0$ was obtained by Y. Tsutsumi in \cite{Ts1} and T. Cazenave, F. Weissler in \cite{cw}.
It is important to note that the critical regularity for the scaling argument is $s=0$.
On the other hand, for the periodic boundary problem, local well-posedness for data in $H^s(\mathbb{T}^2),~s>0$ was established by J. Bourgain \cite{bo1}.
In \cite{bo1}, he used the $X^{s,b}$ argument along with the Strichartz estimate which controls the $L^{p}_{t,x}$ norm of the linear solution in term of the $H^s$ initial data (see also \cite{z}), where this estimate is sharp in the sense that the estimate with $s\le 0$ fails.
For further references, see for instance \cite{bo3}, \cite{ta1} and references therein.

%We summarize next some local results for the initial value problem of (\ref{NLS}).
Among other cases, for the case of two dimensional sphere $S^2$, the initial value problem (\ref{NLS}) is known to be locally well-posed in the Sobolev spaces $H^s(S^2)$ for $s>1/4$.
This was shown by N. Burq, P. G\'erard and N. Tzvetkov \cite{bgt2}\footnote{In \cite{bgt2}, the same problem on a compact Riemann manifold has been discussed.}, where they proved the bilinear Strichartz estimates without loss of derivatives.
On the other hand, in \cite{bgt1}, they also proved that the flow map is not uniformly continuous in $H^s(S^2)$ for $0\le s<1/4$, although the standard scaling argument concerning the initial value problem (\ref{NLS}) on the real line case suggests local well-posedness for $s\ge 0$.
Further analysis on such instability of solutions was obtained by V. Banica in \cite{Ha}.
It is worth noticing that these instability phenomena arise to be illustrate by considering the following weighted homogeneous spherical harmonic polynomials
\begin{eqnarray}\label{hom-shperical}
\psi_k(x_1,x_2,x_3)=k^{1/4-s}(x_1+ix_2)^k,\quad
(x_1,x_2,x_3)\in S^2,~k\in\mathbb{N},
\end{eqnarray}
as initial data.
%This along with the frequency modulation technique implies that the map data-solution is not 
Notice that these weighted homogeneous spherical harmonic polynomials give a concrete counterexample to %the bilinear Strichartz estimate obtained by N. Burq, P. G\'erard and N. Tzvetkov (also to 
the restriction estimate obtained by C. D. Sogge \cite{so1}.

It is natural to ask for the well-posedness result at the regularity $s=1/4$.
In the present paper, we consider some problems at the regularity $s=1/4$.

We first recall the spherical harmonics on the sphere $S^2$.
\begin{notation}
We define $A\lesssim B$ by $A\le cB$ for some absolute constant $c>0$.
Similarly, we define $A\sim B$ if and only if $A\lesssim B\lesssim A$.
We set $\langle a\rangle=(1+|a|^2)^{1/2}$ for $a\in \mathbb{R}$. 

We summarize the properties of eigenvalues and eigenfunctions of the Laplacian $\Delta$ on $S^2$.
Let $Y_k^m~(k\in \mathbb{N}\cup\{0\},~m\in\mathbb{Z},~-k\le m\le k)$ be complex valued spherical harmonics of degree $k$ and order $m$.
Then $Y_k^m$ are eigenfunctions corresponding to eigenvalues $-\mu_k^2$, where $\mu_k=\sqrt{k(k+1)}$.
For a fixed $k$, we denote by ${\cal H}_k$ the $(2k+1)$-dimensional space spanned by the functions $Y_k^m$ with $-k\le m\le k$. 
We denote by $P_k$ the orthogonal projection on eigenspace ${\cal H}_k$\footnote{The space $L^2(S^2)$ is the Hilbert sum $\bigoplus_{k=0}^{\infty}{\cal H}_k$, namely, for every function $f\in L^2(S^2)$, there is a unique Fourier decomposition $f=\sum_{k=0}^{\infty}f_k$ for $f_k\in {\cal H}_k$.}.

Let $s\ge 0$.
Denote by $H^s=H^s(S^2)$ the Sobolev space associated to $(I-\Delta)^{s/2}$ equipped with norm
$$
\|u\|_{H^s}=\left(\sum_{k=0}^{\infty}\langle \mu_k\rangle^{2s}\|P_ku\|_{L^2}^2\right)^{1/2}.
$$
Denote by $B^{s}=B^s(S^2)$ the Besov type modulation space equipped with norm
$$
\|u\|_{B^s}=\sum_{k=0}^{\infty}\langle \mu_k\rangle^{s}\|P_ku\|_{L^2}.
$$
\end{notation}

\begin{remark}
\begin{itemize}
\item[(i)]
We only consider the case when $s=1/4$, since the argument in this paper will provide the proof in the case $s>1/4$.
\item[(ii)]
The definition of reasonable Besov space is $B^{s}_{p,q}$ equipped with norm 
$$
\|u\|_{B^s_{p,q}}=\left\|2^{js}\left\|\sum_{2^j-1\le k<2^j}P_ku\right\|_{L^p(S^2)}\right\|_{\ell_{j}^q}
$$
(see \cite{bl}).
%where $P_ku=\sum_{m=-k}^kP_k^mu$.
There is the embedding inclusion relation between $B^s$ and $B^s_{p,q}$, namely 
$$
B^s \hookrightarrow B_{2,1}^s\hookrightarrow H^s\hookrightarrow B_{2,\infty}^s.
$$
From this point of view, the space $B^s$ is a restricted space of the Besov space.
%The space $B^s$ is a modulation of $B^{s}_{2,1}$ in this sense.
\item[(iii)]
When $s=1/4$, the no-weighted homogeneous spherical harmonic functions $\widetilde{\psi_k}(x_1,x_2,x_3)=(x_1+ix_2)^k,~(x_1,x_2,x_3)\in S^2$ satisfies $\widetilde{\psi_k}=c_kY_k^k$ for some constant $c_k$ and
$$
P_k\widetilde{\psi_k}=\widetilde{\psi_k},\quad \widetilde{\psi_k}\in B^{1/4},$$
$$
\|\widetilde{\psi_k}\|_{B^{1/4}}=\|\widetilde{\psi_k}\|_{B_{2,1}^{1/4}}=\|\widetilde{\psi_k}\|_{H^{1/4}}=\|\widetilde{\psi_k}\|_{B_{2,\infty}^{1/4}}\sim 1.
$$
\end{itemize}
\end{remark}

The main result of this paper is local well-posedness of (\ref{NLS}) in $B^{1/4}$.
More precisely, we show the following theorem.
\begin{theorem}\label{LWP_B}
For every $r>0$, there exist $T=T(r)>0$ and a function space $X_T^{1/4}$ such that given $u_0\in \{\phi \in B^{1/4}\mid \|\phi\|_{B^{1/4}}<r\}$ there exists a unique solution $u$ of (\ref{NLS}) with initial data $u_0$ satisfying $u\in C([0,T];B^{1/4})\cap X^{1/4}_T$.
Moreover the solution map $\{\phi \in B^{1/4}\mid \|\phi\|_{B^{1/4}}< r\}\ni u_0\mapsto u\in C([0,T];B^{1/4})\cap X^{1/4}_T$ is Lipschitz continuous.
% from the neighborhood of $u_0$ in $B^{1/4}$ into $X^{1/4}_T$.
\end{theorem}

\begin{remark}
\begin{itemize}
\item[(i)]
The time of existence $T$ obtained in Theorem \ref{LWP_B} depends on the value of $\|u_0\|_{B^{1/4}}$.
For sequence of initial data $\widetilde{\psi_k}~(k\in\mathbb{N})$, Theorem \ref{LWP_B} guarantees that there exists the sequence of solutions $u_k$ and its existence time $T_k$.
We can say that by Theorem \ref{LWP_B} the lifespan has the lower estimate $\liminf_{k \to\infty} T_k\ge c$ for some constant $c>0$.
\item[(ii)]
The result in Theorem \ref{LWP_B} is essentially of the some kind of conditionally local-in-time well-posedness in $B^{1/4}$.
Our proof relies on the a priori estimates of solutions.
The same sprit in carrying out the proof also works for the unconditionally well-posedness in $H^s$ with $s>3/4$.
However, in this paper we do not discuss at all such a problem.
%More precisely, we obtain the following theorem.
%
%\begin{theorem}\label{LWP_H}
%Let $s>3/4$.
%Then for every $u_0\in H^s$, there exists a unique solution $u\in C([0,T];H^s)$ of (\ref{NLS}) with initial data $u_0$, where $T$ of existence depends only on $\|u_0\|_{H^s}$.
%Moreover the solution map $u_0\to u(t)$ is Lipschitz continuous from the neighborhood of $u_0$ in $H^s$ into $C([0,T];H^s)$.
%\end{theorem}
\end{itemize}
\end{remark}

The well-posedness result for (\ref{NLS}) is still unclear.
But it is interesting to consider the cubic nonlinear Schr\"odinger equation without gauge invariance.
If data are homogeneous spherical harmonic functions, then one can obtain the local well-posedness for data in the Sobolev space $H^s_{hom}$ for $s=1/4$, where
$$
H^s_{hom}=H^s\cap \mathrm{span}\{\psi_k\mid k\ge 0\}.
$$

\begin{theorem}\label{thm-toy-NLS}
The initial value problem associated with the cubic nonlinear Schr\"odinger equation
\begin{eqnarray}\label{toy-NLS}
i\partial_t u+\Delta u=u^3,
\end{eqnarray}
is locally well-posed in $H^{1/4}_{hom}$.
Namely, for every $r>0$, there exist $b>1/2,~T=T(r)>0$ and a function space $X_T^{1/4,b}$ such that given $u_0\in \{\phi \in H_{hom}^{1/4}\mid \|\phi\|_{H_{hom}^{1/4}}< r\}$ there exists a unique solution $u$ of (\ref{toy-NLS}) with initial data $u_0$ satisfying $u\in C([0,T];H_{hom}^{1/4})\cap X^{1/4,b}_T$.
Moreover the solution map $\{\phi \in H_{hom}^{1/4}\mid \|\phi\|_{H_{hom}^{1/4}}< r\}\ni u_0\mapsto u\in C([0,T];H_{hom}^{1/4})\cap X^{1/4,b}_T$ is Lipschitz continuous.
%More precisely, for every $u_0\in H_{hom}^{1/4}$, there exists a unique solution $u\in C([0,T];H_{hom}^{1/4})\cap X_T^{1/4,b}$ of (\ref{toy-NLS})
\end{theorem}

\begin{remark}\label{rem-hom}
If $u_j$ are smooth in $\mathrm{span}\{\psi_k\mid k\ge 0\}$ for $j=1,2,3$, then the product $\prod_{j=1}^3u_j$ is within $\mathrm{span}\{\psi_k\mid k\ge 0\}$.
This algebraic property allows us to establish the existence of an appropriate class of local solution $u(t)\in C([0,T];H_{hom}^{s})$ of (\ref{toy-NLS}).
\end{remark}

This paper is organized as follows.
In section 2, we present the useful function spaces and recall the several estimates for the harmonic projection operators.
In section 3, we take the direct estimate for the norm of $\|u(t)\|_{X^s_T}$ along with the energy method under resonance and has several other properties.
Multilinear estimates provide us the several type of a priori estimates.
In section 4, we carry out the local well-posedness result described in Theorem \ref{LWP_B}.
In section 5, we give the proof of Theorem \ref{thm-toy-NLS}.

\noindent
{\it Acknowledgements.}
I thank Nikolay Tzvetkov for discussion and for pointing out the reference \cite{bgt3}.

\section{Definition of function spaces and linear estimates}
We start out by defining the function spaces.
\begin{definition}
Let $s\ge 0$ and $T>0$.
We define the function space $X_T^s$ as the completion of $C_0^{\infty}(\mathbb{R}\times S^2)$ under the norm:
$$
\|u\|_{X_T^s}=\sum_{k=0}^{\infty}\langle \mu_k\rangle^s\sup_{0\le t\le T}\|P_ku(t)\|_{L^2(S^2)}.
$$
%For $s=1/4$, we shall omit the index $s$.
\end{definition}

We also define the function spaces introduced by Bourgain in \cite{bo1} and Burq-G\'erard-Tzvetkov in \cite{bgt3}.
\begin{definition}
For $s\ge 0$ and $b\in\mathbb{R}$, let $X^{s,b}$ be the completion of $C_0^{\infty}(\mathbb{R};H^s(S^2))$ for the norm
$$
\|f\|_{X^{s,b}}=\left(\sum_{k=0}^{\infty}\|\langle \tau+\mu_k^2\rangle^b\langle\mu_k\rangle^s\widehat{P_ku}(\tau)\|_{L^2(\mathbb{R}\times S^2)}^2\right)^{1/2}.
$$
Observe that if $b>1/2$, $X^{s,b}\hookrightarrow C(\mathbb{R};H^s)$ by usual Sobolev embedding theorem.
For $T>0$, we define the restriction function space
$$
X_T^{s,b}=\{f|_{|t|\le T}\mid f\in X^{s,b}\}.
$$
\end{definition}

The following extension of the $L^p$ estimates due to Sogge \cite{so1} established by Burq, G\'erard and Tzvetkov \cite{bgt3} will be a crucial ingredient in our proof.

\begin{theorem}[Burq, G\'erard and Tzvetkov]\label{bilinear}
There exist $C>0$ such that, for all $k_1,k_2\in \mathbb{N}\cup\{0\}$, for all functions $f,g$ on $S^2$,
$$
\|P_{k_1}fP_{k_2}g\|_{L^2(S^2)}\le C\min\{\langle k_1\rangle,\langle k_2\rangle\}^{1/4}\|f\|_{L^2(S^2)}\|g\|_{L^2(S^2)}.
$$
\end{theorem}
Also we shall use the following $L^p$ estimates established by Sogge \cite[Theorem 4.2]{so1}.
\begin{theorem}[Sogge]\label{lem:sogge}
There exist $C>0$ such that, for all $k\in\mathbb{N}\cup\{0\}$, for all $f$ on $S^2$,
$$
\|P_ku\|_{L^p(S^2)}\le C\langle k\rangle^{\alpha_p}\|u\|_{L^2(S^2)},
$$
where,
\begin{itemize}
\item[(i)]
if $2\le p \le 6$, $\alpha_p=(p-2)/4p$,
\item[(ii)]
if $6\le p\le\infty$, $\alpha_p=2(1/4-1/p)$.
\end{itemize}
\end{theorem}

\section{Nonlinear energy estimates under resonance properties}
To show the local well-posedness in $B^{1/4}$, we use the a priori estimates for solution of (\ref{NLS}).
This can be established as follows.
\begin{theorem}\label{apriori}
Let $u$ be a smooth solution of (\ref{NLS}) on $t\in[0,T]$.
Then for all $0<\delta<1$, we have
$$
\|u\|_{X_T^{1/4}}\lesssim \|u_0\|_{B^{1/4}}+\left(\sqrt{\frac{T}{\delta}}+\delta^{1/4}\right)\|u\|_{X_T^{1/4}}^2+\sqrt{T}\|u\|_{X_T^{1/4}}^3.
$$
\end{theorem}

To construct above estimate, first consider a smooth solution $u$ of (\ref{NLS}) with smooth initial data $u_0$.
%$S(-t)=e^{-it\Delta}$ denote the semigroup to the linear Schr\"odinger equation.
We put $v_k=e^{-it\Delta}P_ku=e^{it\mu_k^2}P_ku$, where we suppressed the $t,x$ dependence and wrote $v_{k_j}=v_{k_j}(t,x)$ for notational simplicity.

Since from the fact that $e^{it\Delta}$ is unitary group in $L^2$, we have $\|v_k\|_{L^2}=\|P_ku\|_{L^2}$.
Then $v_k$ satisfies the equivalent formula to (\ref{NLS}):
\begin{eqnarray}\label{eq_v}
\partial_tv_k=-iP_k\sum_{k_1,k_2,k_3}e^{i\phi(k,k_1,k_2,k_3)t}v_{k_1}\overline{v_{k_2}}v_{k_3},
\end{eqnarray}
where the phase function $\phi(k,k_1,k_2,k_3)$ is defined by
$$
\phi(k,k_1,k_2,k_3)=\mu_k^2-\mu_{k_1}^2+\mu_{k_2}^2-\mu_{k_3}^2.
$$
%and we suppressed the $x,t$ dependence in (\ref{eq_v}).
Additionally we may assume that $k\lesssim\max\{k_1,k_2,k_3\}$.
Otherwise, the right-hand side of (\ref{eq_v}) is zero, since the $v_{k_1}\overline{v_{k_2}}v_{k_3}$ is the polynomial on degree at most $k_1+k_2+k_3$, and it is orthogonal to the polynomial function of degree $k$.

Then it formally satisfies the identity
\begin{eqnarray}
%& & \|P_ku(t)\|_{L^2}^2-\|P_ku_0\|_{L^2}^2=\int_0^t
& & \frac{d}{dt}\|P_ku(t)\|_{L^2}^2
 =  2\mathrm{Re}\int_{S^2}\overline{v_k}\partial_tv_k\,dx\nonumber\\
 &= & 2\mathrm{Im}\int_{S^2}\overline{v_k(t)}\sum_{k_1,k_2,k_3}e^{i\phi(k,k_1,k_2,k_3)t}v_{k_1}\overline{v_{k_2}}v_{k_3}\,dx\label{diff}.
\end{eqnarray}

Let us drop the complex number $i$ and simply use $c$.
We decompose (\ref{diff}) into a sum of two terms that the resonant contribution $(k_1,k_2,k_3)\in\sigma_k$ and that the non-resonant contribution $(k_1,k_2,k_3)\not\in\sigma_k$, where
\begin{eqnarray}\label{res-relation}
\sigma_{k}=\left\{(k_1,k_2,k_3)\in(\mathbb{N}\cup\{0\})^3:|\mu_k-\sqrt{\mu_{k_1}^2-\mu_{k_2}^2+\mu_{k_3}^2}|\le 1~\mbox{or}~\mu_k\le 1/\delta\right\}.
\end{eqnarray}
From the contribution of $\sigma_k$, we know that
\begin{eqnarray}\label{support}
\mbox{if $|k-k'|>2$ and $\min\{\mu_k,\mu_{k'}\}>2/\delta$, then $\sigma_k\cap\sigma_{k'}=\emptyset$},
\end{eqnarray}
and
\begin{eqnarray}\label{support_c}
\mbox{if $(k_1,k_2,k_3)\not\in\sigma_k$, then $\phi(k,k_1,k_2,k_3)\ne 0$}.
\end{eqnarray}

Using (\ref{diff}), one has after integration in the time interval $[0,t]$ that
\begin{eqnarray*}
& & \|P_ku(t)\|_{B^{1/4}}^2-\|P_ku_0\|_{B^{1/4}}^2\nonumber\\
 & = &  c\int_0^t\int_{S^2}\langle\mu_k\rangle^{\frac12}\overline{v_k}\sum_{(k_1,k_2,k_3)\in\sigma_{k}}e^{i\phi(k,k_1,k_2,k_3)t'}v_{k_1}\overline{v_{k_2}}v_{k_3}\,dxdt'\label{res}\\
& & +c\int_0^t\int_{S^2}\langle\mu_k\rangle^{\frac12}\overline{v_k}\sum_{(k_1,k_2,k_3)\not\in\sigma_k}e^{i\phi(k,k_1,k_2,k_3)t'}v_{k_1}\overline{v_{k_2}}v_{k_3}\,dxdt'\label{nonres}\\
& = & I_1+I_2.
\end{eqnarray*}

{\it Contribution of $(k_1,k_2,k_3)\in\sigma_k$.}
It is enough to consider the case where $k_1$ is the maximal modulation of $k_j~(j=1,2,3)$.
By two applications of Theorem \ref{bilinear}, we obtain
\begin{eqnarray*}
|I_1| & \lesssim & \int_0^t\sum_{(k_1,k_2,k_3)\in\sigma_k}\langle \mu_k\rangle^{\frac14}\langle \mu_{k_1}\rangle^{\frac14}\|v_kv_{k_2}\|_{L^2}\|v_{k_1}v_{k_3}\|_{L^2}\,dt'\\
& \lesssim & \int_0^t\|v_k\|_{B^{1/4}}\sum_{(k_1,k_2,k_3)\in\sigma_k}\prod_{j=1}^3\|v_{k_j}\|_{B^{1/4}}\,dt'.
\end{eqnarray*}

{\it Contribution of $(k_1,k_2,k_3)\not\in\sigma_k$.}
Since $\phi(k,k_1,k_2,k_3)\ne 0$ by (\ref{support_c}), we write after integration by parts that% the contribution of $(k_1,k_2,k_3)\not\in\sigma_k$ is
\begin{eqnarray*}
I_2 & = & c\left[\int_{S^2}\frac{\langle\mu_k\rangle^{\frac12}e^{i\phi(k,k_1,k_2,k_3)t'}}{\phi(k,k_1,k_2,k_3)}\overline{v_k}\sum_{(k_1,k_2,k_3)\not\in\sigma_k}v_{k_1}\overline{v_{k_2}}v_{k_3}\right]_{t'=0}^{t'=t}\label{nonreso-1}\\
& & +  c\int_0^t\int_{S^2}\frac{\langle\mu_k\rangle^{\frac12}e^{i\phi(k,k_1,k_2,k_3)t'}}{\phi(k,k_1,k_2,k_3)}\partial_t\left(\overline{v_k}\sum_{(k_1,k_2,k_3)\not\in\sigma_k}v_{k_1}\overline{v_{k_2}}v_{k_3}\right)\,dt'\label{nonreso-2}\\
& = & I_{21}+I_{22}.
\end{eqnarray*}

We consider the contribution of $I_{21}$ first.
Note that from the non-resonance relation (\ref{res-relation}), we have
$$
|\phi(k,k_1,k_2,k_3)|\gtrsim \left\langle \mu_k-\sqrt{\mu_{k_1}^2-\mu_{k_2}^2+\mu_{k_3}^2}\right\rangle\langle\mu_k\rangle
$$
and $\langle\mu_k\rangle\gtrsim 1/\delta$ within the domain of summation.
Then, using Theorem \ref{bilinear}, we have that the contribution of $I_{21}$ is bounded by
\begin{eqnarray*}
|I_{21}|& \lesssim & \delta^{1-\varepsilon}\|P_ku(t)\|_{B^{1/4}}\sum_{(k_1,k_2,k_3)\not\in\sigma_k}\frac{\prod_{j=1}^3\|P_{k_j}u(t)\|_{B^{1/4}}}{\left\langle \mu_k-\sqrt{\mu_{k_1}^2-\mu_{k_2}^2+\mu_{k_3}^2}\right\rangle\langle\mu_k\rangle^{\varepsilon}}\\
& &  +\delta^{1-\varepsilon}\|P_ku_0\|_{B^{1/4}}\sum_{(k_1,k_2,k_3)\not\in\sigma_k}\frac{\prod_{j=1}^3\|P_{k_j}u_0\|_{B^{1/4}}}{\left\langle \mu_k-\sqrt{\mu_{k_1}^2-\mu_{k_2}^2+\mu_{k_3}^2}\right\rangle\langle\mu_k\rangle^{\varepsilon}}.
\end{eqnarray*}

Secondly, we consider the contribution of $I_{22}$.
According to the symmetry, we divide the proof into two cases: the time derivative falls on the factor that $v_k$ and that $v_{k_1}$ .%of $\overline{v_k}v_{k_1}\overline{v_{k_2}}v_{k_3}$.

When the time derivative falls on the factor $v_k$, $\partial_tv_k$ satisfies the equation (\ref{eq_v}).
Apply Theorem \ref{bilinear} as above, then the contribution of this case to $I_{22}$ is bounded by
\begin{eqnarray*}
& &  c\int_0^t\frac{\|P_k(|u|^2u)\|_{B^{1/4}}\sum_{(k_1,k_2,k_3)\not\in\sigma_k}\prod_{j=1}^3\|P_{k_j}u\|_{B^{1/4}}}{\left\langle \mu_k-\sqrt{\mu_{k_1}^2-\mu_{k_2}^2+\mu_{k_3}^2}\right\rangle\langle\mu_k\rangle}\,dt'\\
& \lesssim &\int_0^t\frac{\|P_k(|u|^2u)\|_{L^2}}{\langle\mu_k\rangle^{\frac34-\varepsilon}}\,dt'\sup_{0\le t'\le t}\sum_{(k_1,k_2,k_3)\not\in\sigma_k}\frac{\prod_{j=1}^3\|P_{k_j}u\|_{B^{1/4}}}{\left\langle \mu_k-\sqrt{\mu_{k_1}^2-\mu_{k_2}^2+\mu_{k_3}^2}\right\rangle\langle\mu_k\rangle^{\varepsilon}}.
\end{eqnarray*}

On the other hand, when the time derivative falls on the factor $v_{k_1}$, we divide two cases that $\mu_k\gtrsim \mu_{k_1}$ and that $\mu_k\ll \mu_{k_1}$.

We consider the contribution of $\mu_k\gtrsim \mu_{k_1}$.
By evaluating the weight $\langle \mu_k\rangle^{-1/4}$ by $\langle\mu_{k_1}\rangle^{-1/4}$ in the expression above, the same arguments as in the previous case lead to
$$
\left|\int_{S^2}\frac{\overline{P_ku}P_{k_1}(|u|^2u)\overline{P_{k_2}u}P_{k_3}u}{\langle\mu_k\rangle}\,dx\right|\lesssim \frac{\|P_ku\|_{B^{1/4}}\|P_{k_1}(|u|^2u)\|_{L^2}\|P_{k_2}u\|_{B^{1/4}}\|P_{k_3}u\|_{B^{1/4}}}{\langle\mu_k\rangle^{\varepsilon}\langle\mu_{k_1}\rangle^{\frac34-\varepsilon}}.
$$

Next we consider the contribution of $\mu_k\ll \mu_{k_1}$.
Applying Theorem \ref{bilinear} as above, one has that
$$
\left|\int_{S^2}\frac{\overline{P_ku}P_{k_1}(|u|^2u)\overline{P_{k_2}u}P_{k_3}u}{\left\langle \mu_k-\sqrt{\mu_{k_1}^2-\mu_{k_2}^2+\mu_{k_3}^2}\right\rangle\langle\mu_k\rangle}\,dx\right|\lesssim \frac{\|P_ku\|_{B^{1/4}}\|P_{k_1}(|u|^2u)\|_{L^2}\|P_{k_2}u\|_{B^{1/4}}\|P_{k_3}u\|_{B^{1/4}}}{\left\langle \mu_k-\sqrt{\mu_{k_1}^2-\mu_{k_2}^2+\mu_{k_3}^2}\right\rangle\langle\mu_{k}\rangle^{\frac34}}.
$$

We put together all the estimates above.
% and prove the following estimate.
Taking the square root (which is concave) and using Cauchy-Schwarz inequality, we obtain the estimate
\begin{eqnarray*}
& & \|u\|_{X_T^{1/4}}\lesssim \|u_0\|_{B^{1/4}}\label{e1}\nonumber\\
& & + \left(\sum_{k=0}^{\infty}\|P_ku\|_{L_t([0,T];B^{1/4})}\right)^{1/2}
\left(\int_0^T\sum_{{\scriptstyle k\ge 0}\atop{\scriptstyle (k_1,k_2,k_3)\in\sigma_k}}\prod_{j=1}^3\|P_{k_j}u(t')\|_{B^{1/4}}\,dt'\right)^{1/2}\label{e2}\\
& & +\delta^{1/2-\varepsilon}\left(\sum_{k=0}^{\infty}\|P_ku\|_{L_t([0,T];B^{\frac14})}\right)^{\frac12}
\left(\sum_{{\scriptstyle k\ge 0}\atop{\scriptstyle (k_1,k_2,k_3)\not\in\sigma_k}}\frac{\prod_{j=0}^3\|P_{k_j}u\|_{L_t^{\infty}([0,T];B^{1/4})}}{\left\langle\mu_k-\sqrt{\mu_{k_1}^2-\mu_{k_2}^2+\mu_{k_3}^2}\right\rangle\langle\mu_k\rangle^{\varepsilon}}\right)^{\frac12}\label{e3}\\
& &
+\left(\int_0^T\sum_{k=0}^{\infty}\frac{\|P_k(|u|^2u)\|_{L^2}}{\langle\mu_k\rangle^{\frac34-\varepsilon}}\,dt'\right)^{\frac12}
\left(\sum_{{\scriptstyle k\ge 0}\atop{\scriptstyle (k_1,k_2,k_3)\not\in\sigma_k}}\frac{\prod_{j=1}^3\|P_{k_j}u\|_{L_t([0,T];B^{\frac14})}}{\left\langle \mu_k-\sqrt{\mu_{k_1}^2-\mu_{k_2}^2+\mu_{k_3}^2}\right\rangle\langle\mu_k\rangle^{\varepsilon}}\right)^{\frac12}\label{e4}\\
& &
+\left(\sum_{k=0}^{\infty}\|P_ku\|_{L_t^{\infty}([0,T];B^{\frac14})}\right)^{\frac12}\times\nonumber\\
& & \quad\times\left(\int_0^T \sum_{{\scriptstyle k\gtrsim k_1\ge 0}\atop{\scriptstyle (k_1,k_2,k_3)\not\in\sigma_k}}\frac{\|P_{k_1}(|u|^2u)\|_{L^2}\|P_{k_2}u\|_{B^{\frac14}}\|P_{k_3}u\|_{B^{\frac14}}}{\left\langle\mu_k-\sqrt{\mu_{k_1}^2-\mu_{k_2}^2+\mu_{k_3}^2}\right\rangle\langle\mu_k\rangle^{\varepsilon}\langle\mu_{k_1}\rangle^{\frac34-\varepsilon}}\,dt'\right)^{\frac12}\label{e5}\\
& & +\left(\sum_{k=0}^{\infty}\|P_ku\|_{L_t^{\infty}([0,T];B^{\frac14})}\right)^{\frac12}\times\nonumber\\
& & \quad \times\left(\int_0^T \sum_{{\scriptstyle k_1\gg k\ge 0}\atop{\scriptstyle (k_1,k_2,k_3)\not\in\sigma_k}}\frac{\|P_{k_1}(|u|^2u)\|_{L^2}\|P_{k_2}u\|_{B^{\frac14}}\|P_{k_3}u\|_{B^{\frac14}}}{\left\langle\mu_k-\sqrt{\mu_{k_1}^2-\mu_{k_2}^2+\mu_{k_3}^2}\right\rangle\langle\mu_{k}\rangle^{\frac34}}\,dt'\right)^{\frac12},\label{e6}
\end{eqnarray*}
at least.
Since $\sum_{k=0}^{\infty}\|P_ku\|_{L_t([0,T];B^{1/4})}=\|u\|_{X_T^{1/4}}$, it is sufficient to estimate the following five terms, $J_i,~i=1,\cdots,5$:
$$
J_1=\int_0^T\sum_{{\scriptstyle k\ge 0}\atop{\scriptstyle (k_1,k_2,k_3)\in\sigma_k}}\prod_{j=1}^3\|P_{k_j}u(t')\|_{B^{1/4}}\,dt',
$$
$$
J_2=\sum_{{\scriptstyle k\ge 0}\atop{\scriptstyle (k_1,k_2,k_3)\not\in\sigma_k}}\frac{\prod_{j=1}^3\|P_{k_j}u\|_{L_t^{\infty}([0,T];B^{1/4})}}{\left\langle\mu_k-\sqrt{\mu_{k_1}^2-\mu_{k_2}^2+\mu_{k_3}^2}\right\rangle\langle\mu_k\rangle^{\varepsilon}},
$$
$$
J_3=\int_0^T\sum_{k=0}^{\infty}\frac{\|P_k(|u|^2u)\|_{L^2}}{\langle\mu_k\rangle^{\frac34-\varepsilon}}\,dt',
$$
$$
J_4=\int_0^T \sum_{{\scriptstyle k\gtrsim k_1\ge 0}\atop{\scriptstyle (k_1,k_2,k_3)\not\in\sigma_k}}\frac{\|P_{k_1}(|u|^2u)\|_{L^2}\|P_{k_2}u\|_{B^{\frac14}}\|P_{k_3}u\|_{B^{\frac14}}}{\left\langle\mu_k-\sqrt{\mu_{k_1}^2-\mu_{k_2}^2+\mu_{k_3}^2}\right\rangle\langle\mu_k\rangle^{\varepsilon}\langle\mu_{k_1}\rangle^{\frac34-\varepsilon}}\,dt',
$$
$$
J_5=\int_0^T \sum_{{\scriptstyle k_1\gg k\ge 0}\atop{\scriptstyle (k_1,k_2,k_3)\not\in\sigma_k}}\frac{\|P_{k_1}(|u|^2u)\|_{L^2}\|P_{k_2}u\|_{B^{\frac14}}\|P_{k_3}u\|_{B^{\frac14}}}{\left\langle\mu_k-\sqrt{\mu_{k_1}^2-\mu_{k_2}^2+\mu_{k_3}^2}\right\rangle\langle\mu_{k}\rangle^{\frac34}}\,dt'.
$$

Subsequently, we estimate $J_i$ for each separately.

{\it Estimate of $J_1$.}
For $\mu_k>1/\delta$, the each support $\sigma_k$ is disjoint by (\ref{support}).
Therefore $J_2$ is bounded by 
$$
J_1 \lesssim \frac{T}{\delta}\sum_{k_1,k_2,k_3}\prod_{j=1}^3\|P_{k_j}u\|_{L_t^{\infty}([0,T];B^{\frac14})}\lesssim \frac{T}{\delta}\|u\|_{X_T^{\frac14}}^3.
$$

{\it Estimate of $J_2$.}
Due to $\left\langle\mu_k-\sqrt{\mu_{k_1}^2-\mu_{k_2}^2+\mu_{k_3}^2}\right\rangle^{-1}\langle\mu_k\rangle^{-\varepsilon}\in\ell_{k}^{1}$, we have that $J_2$ is bounded by
$$
J_2\lesssim\sum_{k_1,k_2,k_3}\prod_{j=1}^3\|P_{k_j}u\|_{L_t^{\infty}([0,T];B^{\frac14})}
\lesssim \|u\|_{X_T^{\frac14}}^3.
$$

{\it Estimate of $J_3$.}
%So it is suffices to estimate the first part of (\ref{e4}) and show the following the estimate
%$$
%\sup_{0\le t\le T}\sum_{k=0}^{\infty}\frac{\|P_k(|u|^2u)\|_{L^2}}{\langle\mu_k\rangle^{\frac34-\varepsilon}}\lesssim \|u\|_{X_T^{\frac14}}^3.
%$$
Cauchy-Schwarz inequality provides $J_3\lesssim \|u\|_{L_t([0,T],L_x^6)}^3$.
%$$
%\sum_{k=0}^{\infty}\frac{\|P_k(|u|^2u)\|_{L^2}}{\langle\mu_k\rangle^{\frac34-\varepsilon}}\lesssim \|u\|_{L^6}^3,
%$$
%which is equal to
Then we have that
$$
J_3\lesssim T\left\|\sum_{k=0}^{\infty}P_ku\right\|_{L_t^{\infty}([0,T];L^6)}^3\lesssim \left(\sum_{k=0}^{\infty}\|P_ku\|_{L_t^{\infty}([0,T];L^6)}\right)^3.
$$
Now we use Sobolev inequality in Theorem \ref{lem:sogge} to obtain
$$
J_3\lesssim T\left(\sum_{k=0}^{\infty}\langle\mu_k\rangle^{\frac16}\|P_ku\|_{L_t^{\infty}([0,T];L^2)}\right)^3\lesssim T\|u\|_{X_T^{\frac16}}^3.
$$
%Therefore, we see that
%$$
%|(\ref{e4})|\lesssim T^{\frac12}\|u\|_{X_T^{\frac14}}^3.
%$$

{\it Estimate of $J_4$.}
%The proof is we use the same argument as above.
%The first part of (\ref{e5}) is equal to $\|u\|_{X_T^{\frac14}}^{\frac12}$.
It suffices to show that
\begin{eqnarray*}\label{e51}
\sup_{0\le t\le T}\sum_{{\scriptstyle k\gtrsim k_1\ge 0}\atop{\scriptstyle (k_1,k_2,k_3)\not\in\sigma_k}}\frac{\|P_{k_1}(|u|^2u)\|_{L^2}\|P_{k_2}u\|_{B^{\frac14}}\|P_{k_3}u\|_{B^{\frac14}}}{\left\langle\mu_k-\sqrt{\mu_{k_1}^2-\mu_{k_2}^2+\mu_{k_3}^2}\right\rangle\langle\mu_k\rangle^{\varepsilon}\langle\mu_{k_1}\rangle^{\frac34-\varepsilon}}
\lesssim \|u\|_{X_T^{\frac14}}^5.
\end{eqnarray*}
The sum $\Sigma_{k}$ is uniformly bounded due to $\left\langle\mu_k-\sqrt{\mu_{k_1}^2-\mu_{k_2}^2+\mu_{k_3}^2}\right\rangle^{-1}\langle\mu_k\rangle^{-\varepsilon}\in\ell_k^1$.
Now we apply same argument as in the previous proof for $J_3$.
Therefore
$$
J_4\lesssim T\|u\|_{X_T^{\frac14}}^5.
$$

{\it Estimate of $J_5$.}
%From the previous estimate of $J_4$, it suffices to prove that
%\begin{eqnarray*}\label{e61}
%\sup_{0\le t\le T}\sum_{{\scriptstyle k_1\gg k\ge 0}\atop{\scriptstyle (k_1,k_2,k_3)\not\in\sigma_k}}\frac{\|P_{k_1}(|u|^2u)\|_{L^2}\|P_{k_2}u\|_{B^{\frac14}}\|P_{k_3}u\|_{B^{\frac14}}}{\left\langle\mu_k-\sqrt{\mu_{k_1}^2-\mu_{k_2}^2+\mu_{k_3}^2}\right\rangle\langle\mu_{k}\rangle^{\frac34}}
%\lesssim \|u\|_{X_T^{\frac14}}^5.
%\end{eqnarray*}
%Using Lemma \ref{int}, 
A simple calculation shows that $J_5$ is bounded by
\begin{eqnarray*}\label{e62}
J_5\lesssim T\sup_{0\le t\le T}\sum_{k_1,k_2,k_3}\frac{\|P_{k_1}(|u|^2u)\|_{L^2}\|P_{k_2}u\|_{B^{\frac14}}\|P_{k_3}u\|_{B^{\frac14}}}{\left\langle\mu_{k_1}^2-\mu_{k_2}^2+\mu_{k_3}^2\right\rangle^{\frac38-\varepsilon}}.
\end{eqnarray*}
Now by Cauchy-Schwarz inequality with the fact $\langle \mu_{k_1}^2+\mu_{k_2}^2+\mu_{k_3}^2\rangle^{-3/8+\varepsilon}\in\ell_{k_1}^2$, we see that
$$
J_5 \lesssim T\|u\|_{L_t([0,T];L^6)}^3.
$$
From the previous step in the proof for $J_3$, we have that
$$
J_5\lesssim T\|u\|_{X_T^{\frac14}}^5.
$$
This completes the proof of Theorem \ref{apriori}.
\qed

Repeating a similar computation on the difference of two solutions, we have the following theorem.

\begin{theorem}\label{apriori-diff}
Let $u$ and $v$ be two smooth solutions of (\ref{NLS}) on $t\in[0,T]$ with initial data $u_0$ and $v_0$, respectively.
Then for all $0<\delta<1$, we have
$$
\|u-v\|_{X_T^{\frac14}}\lesssim \|u_0-v_0\|_{B^{\frac14}}+\left(\sqrt{\frac{T}{\delta}}+\delta^{1/4}\right)M\|u-v\|_{X_T^{\frac14}}+\sqrt{T}M^2\|u-v\|_{X_T^{\frac14}}.
$$
where $M=\max\left\{\|u\|_{X_T^{1/4}},\|v\|_{X_T^{1/4}}\right\}$.
\end{theorem}

\section{Proof of Theorem \ref{LWP_B}}
Now we explain how to use Theorems \ref{apriori} and \ref{apriori-diff} to establish local solutions.
%In this section, we put all estimates in previous section and prove Theorem \ref{LWP_B}.

Let $u_{0,n}$ be a sequence in $H^{\infty}(S^2)$ such that $u_{0,n}\to u_0$ in $B^{1/4}$ as $n\to \infty$.
For $u(x,0)=u_{0,n}(x)$ in (\ref{NLS}), we have global smooth solutions of (\ref{NLS}) which we denote by $u_n$.

Let $T>0$ be a positive constant with $T<1$ to be determined later.
First we use the continuity argument to show that $u_n$ satisfies the a priori estimate obtained in Theorem \ref{apriori} on the time interval $[0,T]$ independent of $n$.
Namely, we have
$$
\|u_n\|_{X_T^{1/4}}\le C\|u_{0,n}\|_{B^{1/4}}+C\left(\sqrt{\frac{T}{\delta}}+\delta^{1/4}\right)\|u_n\|_{X_T^{1/4}}^2+C\sqrt{T}\|u_n\|_{X_T^{1/4}}^3.
$$
We choose $R>0$ so large that $\|u_0\|_{B^{1/4}}\le R,~\|u_{0,n}\|_{B^{1/4}}\le R$ for all $n$.
Since $\|u_n\|_{X_t^{1/4}}$ is continuous in $t$, there exist small $\delta>0$ and  a time interval $[0,T']$ such that $\|u_n\|_{X_{T'}^{1/4}}\le 2CR$.
Then we have by Theorem \ref{apriori}
\begin{eqnarray*}
\|u_n\|_{X_{T'}^{1/4}} & \le & C\|u_{0,n}\|_{B^{1/4}}+C\left(\sqrt{\frac{T'}{\delta}}+\delta^{1/4}\right)\|u_n\|_{X_{T'}^{1/4}}^2+C\sqrt{T'}\|u_n\|_{X_{T'}^{1/4}}^3\\
& \le & CR+\frac{CR}{2}=\frac{3}{2}CR.
\end{eqnarray*}
Therefore we can choose $T<1$ so small that
$$
\|u_n\|_{X_T^{1/4}}\le 2R
$$
where $T$ depends only on $C$ and $R$.

Moreover, we repeat a similar argument on difference of two solutions.
Theorem \ref{apriori-diff} yields
$$
\|u_{n}-u_{m}\|_{X_T^{1/4}}\le C'\|u_{0,n}-u_{0,m}\|_{B^{1/4}}
$$
for some $C'>0$.
Hence by the Ascoli-Arzel'a compactness theorem we construct convergent sequences $\{u_{n(k)}\}\subset\{u_n\}$ by taking subsequences.
Thus we obtain a solution $u$ of (\ref{NLS}) satisfying
$$
u\in X_T^{1/4}\subset L^{\infty}([0,T];B^{1/4}),\quad \|u\|_{X_T^{1/4}}\le 2 R,
$$
$$
\lim_{k\to\infty}\|u_{n(k)}-u\|_{X_T^{1/4}}=0.
$$

We now present the proof of uniqueness of solution to (\ref{NLS}).
Let $u$ and $v$ be two solutions of (\ref{NLS}) in $X_T^{1/4}$ with the same initial data $u_0$.
We evaluate the $X_{T}^{1/4}$ norm of $u-v$ in the same way as above by using Theorem \ref{apriori-diff}.
If we choose $\delta>0,~T'>0$ small, then
$$
\|u-v\|_{X_{T'}^{1/4}}\le \eta\|u-v\|_{X_{T'}^{1/4}}
$$
for some $0<\eta<1$, which shows $u=v$ on $[0,T']$.
By repeating this procedure, we obtain $u=v$ on $[0,T]$.

Moreover repeating the proof of Theorem \ref{apriori} with the functional
$$
\|P_ku(t_1)-P_ku(t_2)\|_{L^2}^2=\int_{t_2}^{t_1}\frac{d}{dt'}\|P_ku(t')-P_ku(t_2)\|_{L^2}^2\,dt',
$$
we obtain that
$$
\|u(t_2)-u(t_1)\|_{B^{1/4}}\le C\left(\sqrt{\frac{t_1-t_2}{\delta}}+\delta^{1/4}\right)\|u\|_{X_{T}^{1/4}}^2+C\sqrt{t_1-t_2}\|u\|_{X_T^{1/4}}^3
$$
whenever $0\le t_2\le t_1\le T$ and $0<\delta<1$.
Thus if $t_1\to t_2$ and $\delta \to 0$, that $u(t_2)\to u(t_1)$ in $B^{1/4}$.
Hence $u\in C([0,T];B^{1/4})$.

The continuous dependence of solution on initial data is proven in the same way as in the proof of existence of solution.
\qed

\section{Proof of Theorem \ref{thm-toy-NLS}}
In this section, we will solve the integral equation
$$
u(t)=e^{it\Delta}u_0-i\int_0^te^{i(t-t')\Delta}u^3(t')\,dt',\quad u_0\in H_{hom}^{1/4},
$$
in $C([0,T];H_{hom}^{1/4})\cap X_T^{1/4,b}$ for appropriate $b>1/2$.

\begin{remark}
A similar argument as in \cite{bgt3} provides us the local well-posedness result of the Cauchy problem (\ref{toy-NLS}) with no change to their proof.
More precisely, if $u_0\in H^s$ with $s>1/4$, there exists a unique solution $u(t)\in H^s$ for local in time.
As it was remarked in Remark \ref{rem-hom}, it is easy to see that we can have the local result in the homogeneous space $H^s_{hom}$ for $s>1/4$.
\end{remark}

We summarize the properties of the space $X^{s,b}$.
\begin{lemma}\label{X-space}
Let $s\ge 0,~1/2<b<b'<1$ and $0<T<1$.
There exists $c>0$ such that
\begin{eqnarray*}
\|\theta_T(t) e^{it\Delta}u_0\|_{X^{s,b}}\le cT^{1/2-b}\|u_0\|_{H^s},
\end{eqnarray*}
\begin{eqnarray*}
\left\|\theta(t) \int_0^te^{i(t-t')\Delta}f(t')\,dt'\right\|_{X^{s,b}}\le c\|f\|_{X^{s,b-1}},
\end{eqnarray*}
$$
\|\theta_Tf\|_{X^{s,b}}\le cT^{1/2-b}\|f\|_{X^{s,b}},
$$
and
\begin{eqnarray*}
\|\theta_Tf\|_{X^{s,b-1}}\le cT^{b'-b}\|f\|_{X^{s,b'-1}},
\end{eqnarray*}
where $\theta\in C_0^{\infty}(\mathbb{R}),~0\le \theta\le 1,~\theta=1$ near $0$, $\mathrm{supp}\theta\subset [-1,1]$, and $\theta_T(t)=\theta(T^{-1}t)$.
\end{lemma}
The proof of Lemma \ref{X-space} is given by Bourgain in \cite{bo1}, and Burq-G\'erard-Tzvetkov in \cite{bgt3}.

We need to have an estimate that takes the nonlinearity $u^3$ in $X^{s,b-1}$.
We comment on the bilinear Strichartz estimate constructed in \cite{bgt3}, and prove Theorem \ref{thm-toy-NLS}.
\begin{lemma}\label{bilinear-St}
Let $b>3/8$ and $f_j\in X^{1/4,b}$ be homogeneous spherical functions of the form
\begin{eqnarray}\label{sp}
f_j(t,x)=\sum_{k\ge 0}a_{j,k}(t)\sin^k\theta e^{ik\phi},
\end{eqnarray}
for $j=1,2$, where $\theta,\phi$ are spherical coordinates $x=(\sin\theta\cos\phi,\sin\theta\sin\phi,\cos\theta),~0\le\theta<\pi,~0\le \phi<2\pi$.
Then the following bilinear estimate holds:
$$
\|f_1f_2\|_{L_{t,x}^2}\lesssim \|f_1\|_{X^{0,b}}\|f_2\|_{X^{1/4,b}}.
$$
\end{lemma}

For the proof of Lemma \ref{bilinear-St}, we will need some elementary inequalities.
\begin{lemma}\label{int}
Let $0<\alpha\le\beta$ such that $\alpha+\beta>1$ and $\varepsilon>0$.
Then
$$
\int_{-\infty}^{\infty}\frac{dt}{\langle t-a\rangle^{\alpha}\langle t-\beta\rangle^{\beta}}\lesssim \frac{1}{\langle \alpha-\beta\rangle^{\gamma}},
$$
where
$$
\gamma=\left\{
\begin{array}{cc}
\alpha+\beta-1, & \mbox{if $\beta<1$},\\
\alpha-\varepsilon, & \mbox{if $\beta=1$},\\
\alpha, & \mbox{if $\beta>1$}.
\end{array}
\right.
$$
\end{lemma}
\noindent
For the proof of Lemma \ref{int}, see \cite{kpv1} for instance.

\noindent
{\it Proof of Lemma \ref{bilinear-St}.}
The proof follows the same general outline as arguments by Bourgain \cite{bo1} and C. E. Kenig, G. Ponce and L. Vega \cite{kpv1}.

%We divide $f,g$ as follows:
%$$
%f=f_++f_-,\quad g=g_++g_-,
%$$
%where
%$$
%f_{\pm}(t,x)=\sum_{n\ge 0}\int_{\mathbb{R}}f_n^{\pm}(t)\sin^n\theta e^{\pm in\phi-i\mu_n^2t},\quad 
%g_{\pm}(t,x)=\sum_{n\ge 0}\int_{\mathbb{R}}g_n^{\pm}(t)\sin^n\theta e^{\pm in\phi-i\mu_n^2t}.
%$$
%Since $\|fg\|_{L^2}=\|f\overline{g}\|_{L^2}$, we have only to estimate the $L^2$ norm of term $f_+g_+$.
%Therefore we may suppose $f_-=g_-=0$.

By the Parseval identity in time variable, we have
\begin{eqnarray*}
\|f_1f_2\|_{L_{t,x}^2} = 
\left\|\int_{\mathbb{R}}\sum_{k_1,k_2\ge 0}\widehat{a_{1,k_1}}(\tau-\tau_1)\widehat{a_{2,k_2}}(\tau_1)\sin^{k_1+k_2+1/2}\theta e^{i(k_1+k_2)\phi}\, d\tau_1 \right\|_{L_{\tau,\theta,\phi}^2}.
\end{eqnarray*}
Using the similar computation to \cite{bo1}, we have that $\|f_1f_2\|_{L^2}^2$ is equal to
\begin{eqnarray*}
& &  \sum_{k_1,k_1',k_2,k_2'}\int_{\mathbb{R}^3}\int_0^{2\pi}\int_0^{\pi}\widehat{a_{1,k_1}}(\tau-\tau_1)\widehat{a_{2,k_2}}(\tau_1)
\overline{\widehat{a_{1,k_1'}}(\tau-\tau_1')\widehat{a_{2,k_2'}}(\tau_1')}\\
& & \times \sin^{k_1+k_1'+k_2+k_2'+1}\theta e^{i(k_1+k_2-k_1'-k_2')\phi}\,d\tau_1 d\tau_1'd\tau d\theta d\phi.
\end{eqnarray*}
It is easy to see that
$$
\int_0^{2\pi}e^{i(k_1+k_2-k_1'-k_2')\phi}\,d\phi=2\pi \delta(k_1+k_2-k_1'-k_2'),
$$
and
$$
\int_0^{\pi}\sin^{k_1+k_1'+k_2+k_2'+1}\theta\,d\theta\sim\frac{1}{\sqrt{k_1+k_1'+k_2+k_2'+1}}.
$$
We write $m=k_1+k_2=k_1'+k_2'$ and so we can conclude that
$$
\|f_1f_2\|_{L^2}^2\lesssim 
\sum_m\frac{1}{\sqrt{2m+1}}\int_{\mathbb{R}}\left(\sum_{k_1}\int_{\mathbb{R}}|\widehat{a_{1,k_1}}(\tau-\tau_1)\widehat{a_{2,m-k_1}}(\tau_1)|\,d\tau_1\right)^2\,d\tau.
$$
Using Cauchy-Schwarz, the inner integral can be estimated as follows:
\begin{eqnarray*}
& & \sum_{k_1}\int_{\mathbb{R}}|\widehat{a_{1,k_1}}(\tau-\tau_1)\widehat{a_{2,m-k_1}}(\tau_1)|\,d\tau_1\\
& \lesssim &  
M\left(\sum_{k_1}\int_{\mathbb{R}}|\widehat{a_{1,k_1}}(\tau-\tau_1)|^2\langle\tau-\tau_1-\mu_{k_1}^2\rangle^{2b} |\widehat{a_{2,m-k_1}}(\tau_1)|^2\langle \tau_1-\mu_{m-k_1}^2\rangle^{2b}\,d\tau_1\right)^{1/2},
\end{eqnarray*}
where
$$
M=\sup_{\tau\in\mathbb{R},m\ge 0}\left(\sum_{k_1}\int_{\mathbb{R}}\langle\tau-\tau_1-\mu_{k_1}^2\rangle^{-2b}\langle \tau_1-\mu_{m-k_1}^2\rangle^{-2b}\,d\tau_1\right)^{1/2}.
$$
A simple computation by Lemma \ref{int} yields
$$
M^2\lesssim \sup_{\tau,m}\sum_{k_1}\langle\tau-\mu_{k_1}^2-\mu_{m-k_1}^2\rangle^{-4b+1}.
$$
For each $\tau$ and $m$, we obtain the two roots $k_1=\alpha_{\pm}(\tau,m)$ of the quadratic equation with respect to $k_1$:
$$
\tau-\mu_{k_1}^2-\mu_{m-k_1}^2=0.
$$
Then we decompose the sum of $\sum_{k_1}$ into two cases:
\begin{itemize}
\item
$|k_1-\alpha_{+}(\tau,m)|<1$ or $|k_1-\alpha_{-}(\tau,m)|<1$,
\item
$|k_1-\alpha_{\pm}(\tau,m)|\ge 1$.
\end{itemize}
Therefore we conclude that
$$
M^2\lesssim 1+\sup_{\tau,m}\sum_{|k_1-\alpha_{\pm}(\tau,m)|\ge 1}\langle k_1-\alpha_{+}(\tau,m)\rangle^{-4b+1}\langle k_1-\alpha_{-}(\tau,m)\rangle^{-4b+1}\lesssim 1,
$$
since by $b>3/8$.

Then it follows that $\|f_1f_2\|_{L^2}^2$ is bounded by
\begin{eqnarray*}
 & &   c\sum_{m,k_1}\frac{1}{\sqrt{2m+1}}\int_{\mathbb{R}^2}|\widehat{a_{1,k_1}}(\tau-\tau_1)|^2\langle\tau-\tau_1-\mu_{k_1}^2\rangle^{2b} |\widehat{a_{2,m-k_1}}(\tau_1)|^2\langle \tau_1-\mu_{m-k_1}^2\rangle^{2b}\,d\tau_1d\tau\\
& & \lesssim \sum_{k_1}\int_{\mathbb{R}}\frac{|\widehat{a_{1,k_1}}(\tau_1)\langle\tau_1-\mu_{k_1}^2\rangle^{b}|^2}{\langle \mu_{k_1}\rangle^{1/4}}\,d\tau_1\sum_{k_2}\int_{\mathbb{R}}|\widehat{a_{2,k_2}}(\tau_2)\langle\tau_2-\mu_{k_2}^2\rangle^{b}|^2\,d\tau_2\\
& & =c\|f_1\|_{X^{0,b}}^2\|f_2\|_{X^{1/4,b}}^2,
\end{eqnarray*}
which completes the proof.
\qed

Lemma \ref{bilinear-St} is equivalent to proving the following result.
\begin{lemma}\label{tri-es}
Let $f_j\in X^{1/4,1/2+}~(1\le j\le 3)$ be of form in (\ref{sp}).
Then for all $b'<5/8,~b>3/8$, the trilinear estimate
\begin{eqnarray}
\left\|\prod_{j=1}^3f_j\right\|_{X^{1/4,b'-1}}\lesssim \prod_{j=1}^3\|f_j\|_{X^{1/4,b}}
\end{eqnarray}
holds.
\end{lemma}
\noindent
{\it Proof.}
First notice that by Remark 1.3 the product function $\prod_{j=1}^3f_j$ preserves of form in (\ref{sp}).
Indeed, one has that there exist $a_{k}(t)~(k\ge 0)$ such that
$$
\prod_{j=1}^3f_j(t,x)=\sum_{k\ge 0}a_{k}(t)\sin^k\theta e^{ik\phi}.
$$
%\begin{eqnarray*}
%& & \widehat{P_k\prod_{j=1}^3f_j}(\tau,x)= \sin^{k}\theta e^{ik\phi}\\
%& & \times\sum_{k_2,k_3}\int_{\mathbb{R}^2}\widehat{a_{1,k-k_2-k_3}}(\tau-\tau_2-\tau_3-\mu_{k-k_2-k_3}^2)\widehat{a_{2,k_2}}(\tau_2-\mu_{k_2}^2)\widehat{a_{3,k_3}}(\tau_3-\mu_{k_3}^2)\,d\tau_2d\tau_3.
%\end{eqnarray*}
Using duality technique, orthogonality argument, Leibniz rule for fractional derivatives and symmetry argument, %a similar argument as the proof of Lemma \ref{bilinear-St}, 
it suffices to show that
\begin{eqnarray}\label{4t}
& & %\int_{\mathbb{R}^3}\sum_{k_2,k_3,k}
%|\widehat{a_{1,k-k_2-k_3}}(\tau-\tau_2-\tau_3-\mu_{k-k_2-k_3}^2)\widehat{a_{2,k_2}}(\tau_2-\mu_{k_2}^2)\widehat{a_{3,k_3}}(\tau_3-\mu_{k_3}^2)\widehat{\overline{a_{4,k}}}(\tau-\mu_k^2)|\,d\tau_2d\tau_3d\tau\\
%& & \lesssim \int_{\mathbb{R}}|\widehat{a_{4,k_4}}(\tau_4+\mu_{k_4}^2)\langle\tau_4-\mu_{k_4}^2\rangle^{1-b'}|^2\,d\tau_j\prod_{j=1}^3\sum_{k_j}\int_{\mathbb{R}}|\widehat{a_{j,k_j}}(\tau_j-\mu_{k_j}^2)\langle\tau_j-\mu_{k_j}^2\rangle^{b}|^2\,d\tau_j.
\left|\int_{\mathbb{R}\times S^2} f_4(t,x)\prod_{j=1}^3f_j(t,x)\,dt dx\right|
\lesssim \|f_1\|_{X^{1/4,b}}\|f_2\|_{X^{1/4,b}}\|f_3\|_{X^{0,b}}\|f_4\|_{X^{0,1-b'}},
\end{eqnarray}
for any $f_4\in X^{0,1-b'}$ of form
$$
f_4(t,x)=\sum_{k\ge 0}a_{4,k}(t)\sin^k\theta e^{ik\phi}.
$$
By two applications of Lemma \ref{bilinear-St}, we have that the left-hand side of (\ref{4t}) is bounded by
$$
\|f_1f_3\|_{L_{t,x}^2}\|f_2f_4\|_{L_{t,x}^2}\lesssim \|f_1\|_{X^{1/4,b}}\|f_3\|_{X^{0,b}}\|f_2\|_{X^{1/4,b}}\|f_4\|_{X^{0,1-b'}},
$$
provided $b>3/8$ and $b'<5/8$.
%where we shall assume non-negative functions $\widehat{a_{j,k_j}},~$ to simplify the exposition.
%We denote respectively the left-hand side by $K$ and the right-hand side by $L$.
%
%Writing $m=k_2+k_3,~\zeta=\tau_2+\tau_3$ and the Cauchy-Schwarz inequality tell us that
%\begin{eqnarray*}
%K& \lesssim &  \sum_m\int_{\mathbb{R}}\left(\sum_{k}\int_{\mathbb{R}}\widehat{a_{1,k-m}}(\tau-\zeta-\mu_{k-m}^2)\widehat{a_{4,k}}(\tau-\mu_k^2)\,d\tau\right)\\
%& & \times
%\left(\sum_{k_2}\int_{\mathbb{R}}\widehat{a_{2,k_2}}(\tau_2-\mu_{k_2}^2)\widehat{a_{3,m-k_2}}(\zeta-\tau_2-\mu_{m-k_2}^2)\,d\tau_2\right)\,d\zeta\\
%& \lesssim & M L,
%\end{eqnarray*}
%where
%$$
%M=\sup_{\zeta,m}\left(\sum_{k,k_2}\int_{\mathbb{R}^2}\frac{d\tau d\tau_2}{\langle \tau-\zeta-\mu_{k-m}^2\rangle^{2b}\langle \tau+\mu_k^2\rangle^{2(1-b')}\langle\tau_2-\mu_{k_2}^2\rangle^{2b}\langle\zeta-\tau_2-\mu_{m-k_2}^2\rangle^{2b}}\right)^{1/2}.
%$$
%A similar argument as the proof of Lemma \ref{bilinear-St} yields
%$$
%M\lesssim \sum_{\zeta,m}\left(\sum_{k,k_2}\langle\zeta+\mu_{k-m}^2+\mu_k^2\rangle^{2(b'-b)-1}\langle\zeta-\mu_{k_2}^2-\mu_{m-k_2}^2\rangle^{-4b+1}\right)^{1/2}\le c.
%$$
Thus lemma follows.
\qed

\noindent
{\it Proof of Theorem \ref{thm-toy-NLS}.}
For positive constants $T$ and $r$, we define
$$
E(T,r)=\left\{u\in C_0(\mathbb{R};H_{hom}^{1/4})\cap X^{s,b}\mid \|\theta_Tu\|_{X^{1/4,b}}\le rT^{1/2-b}\right\}.
$$
We shall show that for appropriate values of $T>0$ and $b>1/2$, the following map $u\to \Phi(u)$ defines a contraction map on $E(T,r)$ for appropriate values of $T$ and $r$:
\begin{eqnarray}\label{u^3-int}
\Phi(u)(t)=e^{it\Delta}u_0-i\int_0^t e^{i(t-t')\Delta}[\theta_T(\theta_Tu^3)](t')\,dt'.
\end{eqnarray}

Using Lemmas \ref{X-space} and \ref{tri-es} combined with the integral formula in (\ref{u^3-int}), we obtain
\begin{eqnarray*}
\|\theta_T\Phi(u)\|_{X^{1/4,b}} & \le & cT^{1/2-b}\|u\|_{H^{1/4}}+cT^{1/2-b}\|\theta_T(\theta_Tu)^3\|_{X^{1/4,b-1}}\\
& \le & cT^{1/2-b}\|u_0\|_{H^{1/4}}+cT^{b'-b+1/2-b}\|(\theta_Tu)^3\|_{X^{1/4,b'-1}}\\
& \le & cT^{1/2-b}\|u_0\|_{H^{1/4}}+cT^{b'-b+1/2-b}\|\theta_Tu\|_{X^{1/4,b}}^3\\
& \le & cT^{1/2-b}\|u_0\|_{H^{1/4}}+cT^{b'-b+3(1/2-b)}r^3
\end{eqnarray*}
for $1/2<b<b'<5/8$.
Set $r=2c\|u_0\|_{H^{1/4}}$.
Choosing $T>0$ small and $b'-b+3(1/2-b)>0$  such that $T^{b'-b+2(1/2-b)}r^2\le 1/2$, we have that
$$
\|\theta_T\Phi(u)\|_{X^{1/4,b}} \le rT^{1/2-b}.
$$
A similar argument to above shows that
$$
\|\theta_T\Phi(u_1)-\theta_T\Phi(u_2)\|_{X^{1/4,b}}\le \frac{1}{2}\|\theta_Tu_1-\theta_Tu_2\|_{X^{1/4,b}}.
$$
Notice that $X^{s,b} \hookrightarrow C(\mathbb{R};H^s)$ when $b>1/2$.
Therefore, a standard iteration argument shows that we obtain a unique solution $u$ of (\ref{toy-NLS}) satisfying
$$
\theta_Tu\in ([0,T]; H_{hom}^s) \cap X^{1/4,b},~\|\theta_Tu\|_{X^{1/4,b}}\le rT^{1/2-b}.
$$
This proves Theorem \ref{thm-toy-NLS}
\qed

%\begin{eqnarray}\label{toy}
%\left\{
%\begin{array}{l}
%i\partial_tu+\Delta u=P_{sp}(|u|^2u),\\
%u(0,x)=u_0(x)=P_{sp}u_0(x),
%\end{array}
%\right.
%\end{eqnarray}
%where $P_{sp}$ is the projection operator on $\mathrm{span}\{\sin^n \theta e^{\pm in\phi}\mid n\in\mathbb{Z}_+\}$.
%
%We have the following result.
%\begin{theorem}
%For every $u_0\in H_{sp}^{1/4}$, there exists a unique solution $u\in C(0,T]);H_{sp}^{1/4})\cap X^{1/4,b}_T$ of (\ref{toy}) with initial data $u_0$m, where the time $T$ of existence depends only on $\|u_0\|_{H_{sp}^{1/4}}$.
%Moreover the solution map $u_0\to u(t)$ is Lipchitz continuous from the neighborhood of $u_0$ in $H_{sp}^{1/4}$ into $X^{1/4,b}_T$.
%\end{theorem}
%\begin{remark}
%Hence, we are in the remarkable situation.
%By nature of construction of solutions, we have the solution described by a linear combination of homogeneous spherical harmonic functions and its dual functions for $s=1/4$.
%It is presently unknown, as in \cite{}, whether the well-posedness is true for general initial data in $H^{1/4}$.
%\end{remark}
%
%The proof is an adaptation of 

\end{document}